\documentclass[11pt,a4paper]{article}
\usepackage{theorem}
\usepackage{amssymb}

\makeatletter
\long\def\@makecaption#1#2{%
  \vskip\abovecaptionskip
  \sbox\@tempboxa{#1 #2}%
  \ifdim \wd\@tempboxa >\hsize
    #1 #2\par
  \else
    \global \@minipagefalse
    \hb@xt@\hsize{\hfil\box\@tempboxa\hfil}%
  \fi
  \vskip\belowcaptionskip}
\makeatother

\newcommand{\C}{\mathbb{C}}
\newcommand{\Q}{\mathbb{Q}}
\newcommand{\mQ}{\mathbb{Q}}

\newcommand{\Z}{\mathbb{Z}}
\newcommand{\Ql}{\mQ_{\ell}}

\newcommand{\Qbar}{{\overline{\mQ}}}

\newcommand{\cB}{\mathcal{B}}
\newcommand{\cC}{\mathcal{C}}
\newcommand{\calD}{\mathcal{D}}
\newcommand{\cS}{\mathcal{S}}

\newcommand{\Gl}{\mathrm{GL}}
\newcommand{\Sl}{\mathrm{SL}}
\newcommand{\SP}{\mathrm{Sp}}

\newcommand{\Alt}{\mathrm{Alt}}

\newcommand{\determ}{\mathrm{det}}

\newcommand{\divis}{\mathrm{div}}

\newcommand{\End}{\mathrm{End}}
\newcommand{\Eo}{\mathrm{End}^{0}}
\newcommand{\et}{\mathrm{\acute{e}t}}

\newcommand{\Hdg}{\mathrm{Hg}}
\newcommand{\HH}{\mathrm{H}}
\newcommand{\Hom}{\mathrm{Hom}}
\newcommand{\Id}{\mathrm{Id}}

\newcommand{\Matr}{\mathrm{M}}

\newcommand{\Norm}{\mathrm{Nm}}
\newcommand{\Nrd}{\mathrm{Nrd}}
\newcommand{\OO}{\mathrm{O}}

\newcommand{\Res}{\mathrm{Res}}
\newcommand{\Stand}{\mathrm{St}}
\newcommand{\Tr}{\mathrm{Tr}}
\newcommand{\Uu}{\mathrm{U}}

\newcommand{\Gdiv}{G_\divis}

\newcommand{\punt}{{\raisebox{.25ex}{\textbf{.}}}}
\newcommand{\Qed}{\unskip\nolinebreak\hfill\hbox{\quad $\square$}}
\newcommand{\mijnpijl}{\relbar\joinrel\kern -.25em\rightarrow}

\newcommand{\mijnlangpijl}{\relbar\joinrel\longrightarrow}

\newcommand{\verylongarrow}{\relbar\joinrel\relbar\joinrel\mijnlangpijl}
\newcommand{\isomarrow}{\stackrel{\raise.25ex\hbox{$\sim$}}{\smash{\mijnpijl}}}

\newcommand{\hooklongrightarrow}{\lhook\joinrel\longrightarrow}

\newcounter{numero}
\newcommand{\hrom}[1]{\setcounter{numero}{#1}\Roman{numero}}

\begin{document}
\begin{center}
{\bfseries\Large{Weil classes on abelian varieties}}\\[5pt]
{\large B.J.J.~Moonen and Yu.G.~Zarhin\footnote{Supported by the National Science Foundation}}
\end{center}
\vspace{.5\baselineskip plus .5\baselineskip minus .5\baselineskip}

\noindent\refstepcounter{section}\label{WCsetup}%
\textbf{\thesection.~Introduction.}\quad Let $X$ be a smooth projective variety over $\C$.
We write $\cB^\punt(X) = \oplus_i \cB^i(X)$ for its Hodge ring (so $\cB^i(X)
\subseteq \HH^{2i}(X,\Q)$ is the subspace of Hodge classes). We call a class $c
\in \cB^\punt(X)$ {\sl decomposable\/} if it lies in the subalgebra
$\calD^\punt(X) \subseteq \cB^\punt(X)$ generated by divisor classes. The 
non-decomposable Hodge classes are called \textsl{exceptional\/} classes. As a
consequence of the Lefschetz theorem on $(1,1)$ classes, the decomposable
classes are algebraic; in particular, if $\cB^\punt(X) = \calD^\punt(X)$  then
the Hodge conjecture for $X$ is true. Thus one is  naturally led to the
question whether there exist varieties $X$ for which $\calD^\punt(X) \subsetneq
\cB^\punt(X)$.

That this question has a positive answer was shown in \cite{Pohl}, where an
example (due to Mumford) was given of an abelian fourfold $X$ of CM-type for
which $\calD^2(X) \subsetneq \cB^2(X)$. A few years later,  Weil \cite{Weil} 
showed that the essential ingredient in Mumford's example is
the fact that $\Eo(X)$ contains an imaginary quadratic subfield $k$ which acts
with multiplicities $n_\sigma = n_{\sigma^\prime} = 2$ on the tangent space of
$X$. To the action of this field $k$ one associates a 2-dimensional subspace
$W_k \subseteq \HH^4(X,\Q)$ which, grace to the assumption $n_\sigma =
n_{\sigma^\prime}$, consists of Hodge classes. Moreover, Weil showed that for
``generic'' abelian fourfolds with complex multiplication by $k$ (subject to
the condition $n_\sigma = n_{\sigma^\prime}$), the space $W_k$ consists of
exceptional Hodge classes. Later, Shioda \cite{Shioda} constructed an exceptional 
2-dimensional algebraic Weil class on a non-simple abelian fourfold of Fermat type.

In this paper, we want to study to what extent Weil's method to construct
exceptional Hodge classes can be generalized. To make this more explicit,
consider an abelian variety $X$ of dimension $g \geq 1$, and suppose $F$ is a
subfield of $\Eo(X)$, with $1 \in F$ acting as the identity on $X$. Write $V_X
= \HH^1(X,\Q)$ and let $r = 2g/[F:\Q]$. The 1-dimensional $F$-vector space
$$W_F = W_F(X) := \bigwedge^r_F V_X$$
can be identified in a natural manner with a subspace of $\HH^r(X,\Q)$, see
Lemma~\ref{eltstrivonWF}~(\romannumeral1). We call $W_F$ the space of Weil
classes with respect to $F$.

The main questions that we are interested in are
\begin{itemize}
\item[{\bf Q1:}] under what conditions on $F$ does $W_F$ contain, or even
consist of, Hodge classes?
\item[{\bf Q2:}] if $W_F$ contains Hodge classes, under what conditions on $F$
are these exceptional?
\end{itemize}

Similar questions can be asked for Tate classes (in case $X$ is defined over a
number field). We postpone this to section~\ref{Q1Q2Tate}.
\vspace{\baselineskip}

\noindent\refstepcounter{section}\label{allornothing}%
\textbf{\thesection.}\quad The multiplicative group $F^\ast$ acts on 
$\oplus_{i} \HH^i(X,\Q) \supset \cB^\punt(X) \supseteq \calD^\punt(X)$ and on 
the subspace $W_F \subseteq \HH^r(X,\Q)$. If the latter action is given by 
$\rho \colon F^\ast \rightarrow \Gl(W_F)$ then its relation to the natural 
structure of $F$-vector space on
$W_F$ is given by $\rho(f)(w) = f^r \cdot w$ for all $f \in F^\ast$, $w \in W$.
Since $\dim_F(W_F) = 1$, it readily follows that {\sl either all elements of
$W_F$ are Hodge classes, or $0 \in W_F$ is the only Hodge class}, and in the
first case, {\sl either $W_F \setminus \{0\}$ consists entirely of exceptional
classes, or none of the classes in $W_F$ is exceptional}.
\vspace{\baselineskip}

\noindent\refstepcounter{section}\label{HodgedecompWF}%
\textbf{\thesection.}\quad Write $\Sigma_F$ for the set of embeddings 
$F \rightarrow \C$, and for $\sigma \in \Sigma_F$, let $\sigma^\prime$ denote 
its complex conjugate. The action of
$F$ on $V_X$ gives a decomposition of $V_\C = V_X \otimes_\Q \C$ as

$$V_\C = \bigoplus_{\sigma \in \Sigma_F} V_{\C,\sigma} = \bigoplus_{\sigma \in
\Sigma_F} \big(V_{\C,\sigma}^{1,0} \oplus V_{\C,\sigma}^{0,1}\big)\, .$$
The dimension $n_\sigma$ of $V_{\C,\sigma}^{1,0}$ is called the multiplicity of
$\sigma$ on the tangent space of $X$; we have $n_\sigma + n_{\sigma^\prime} =
2g/[F:\Q]$ for all $\sigma \in \Sigma_F$. Similar to the computation in the
proof of \cite[Prop.~4.4]{Del2} (see also \cite[Lemma~2.8]{MZ}), we can write
$$W_F \otimes \C = \big(\bigwedge^r_F V_X \big) \otimes \C = \bigwedge_{F
\otimes_\Q \C}^r V_\C = \bigoplus_{\sigma \in \Sigma_F} \big(\bigwedge_\C^r
V_{\C,\sigma} \big) = \bigoplus_{\sigma \in \Sigma_F}
\big(\bigwedge^{n_\sigma}_\C V_{\C,\sigma}^{1,0} \otimes
\bigwedge^{n_{\sigma^\prime}} V_{\C,\sigma}^{0,1} \big)\, .$$

In view of the remarks in section~\ref{allornothing} we obtain the following answer 
to question Q1.
\vspace{\baselineskip}

\noindent\refstepcounter{section}\label{crit1}%
\textbf{\thesection.~Criterion.}\quad \textsl{If $n_\sigma = n_{\sigma^\prime}$ 
for all $\sigma \in \Sigma_F$ then $W_F$
consists entirely of Hodge classes; if $n_\sigma \neq n_{\sigma^\prime}$ for
some $\sigma \in \Sigma_F$ then the zero class is the only Hodge class in
$W_F$.}
\vspace{\baselineskip}

\noindent\refstepcounter{section}%
\textbf{\thesection.~Remark.}\quad In the case of an imaginary quadratic field $F$ 
this criterion is essentially due to Weil \cite{Weil}; see also Ribet \cite{Ribet}. 
The case of arbitrary CM-fields $F$ was also treated in \cite{SZ1}.
\vspace{\baselineskip}

\noindent\refstepcounter{section}%
\textbf{\thesection.~Remark.}\quad If all simple factors of $X$ are of type \hrom{1}, 
\hrom{2} or \hrom{3} in the Albert classification, then every subfield $F \subseteq
\Eo(X)$ satisfies the condition that $n_\sigma = n_{\sigma^\prime}$ for all
$\sigma \in \Sigma_F$. We can see this as follows: we have an inclusion
$\Hdg(X) \subset \Gl_F(V_X)$ and $\Hdg(X)$ acts on $W_F$ through the $F$-linear
determinant $\determ_F \colon \Gl_F(V_X) \rightarrow \Res_{F/\Q}(\mathbb{G}_{m,F})$. 
If $X$ has no factors of type \hrom{4} then the Hodge group
$\Hdg(X)$ is semi-simple (see \cite[Cor.\ 1.2.2]{Chi}), hence contained in
$\Sl_F(V_X)$, which means that $W_F$ consists of Hodge classes.
\vspace{\baselineskip}

\noindent\refstepcounter{section}\label{decompuptoisog}%
\textbf{\thesection.}\quad Up to isogeny we can decompose $X$ as
\begin{equation} X \sim Y_1^{m_1} \times \cdots \times Y_k^{m_k}\,
,\label{eq:Xdecomp}\end{equation}
where $Y_1, \ldots , Y_k$ ($k \in \Z_{\geq 1}$) are simple, mutually
non-isogenous, abelian varieties, and $m_1, \ldots, m_k \in \Z_{\geq 1}$. The
condition that $1 \in F$ acts as the identity implies that $F$ acts on each
factor $Y_i^{m_i}$, and that the action of $F$ on $X$ is the ``diagonal
action'' w.r.t.\ the decomposition (\ref{eq:Xdecomp}).

Write $r_i = 2 m_i \dim(Y_i) / [F:\Q]$, so that $r = r_1 + \cdots + r_k$. For
each factor $Y_i^{m_i}$ we obtain a 1-dimensional $F$-subspace $W_F(Y_i^{m_i})
\subset \HH^{r_i}(Y_i^{m_i},\Q)$. We claim that the space $W_F(X)$ can be
identified with the tensor product $W_F(Y_1^{m_1}) \otimes_F \ldots \otimes_F
W_F(Y_k^{m_k})$, considered as a subspace (sic) of the K\"unneth component
$$\HH^{r_1}(Y_1^{m_1},\Q) \otimes_\Q \ldots \otimes_\Q \HH^{r_k}(Y_k^{m_k},\Q)
\subset \HH^r(X,\Q)\, .$$
The verification of this statement, which we leave to the reader, is a matter
of linear algebra. Some of our identifications may seem a little unnatural;
this is caused by our desire to view $W_F(X)$ as a {\sl subspace\/} of
$\HH^r(X,\Q)$, rather than a {\sl quotient\/} (cf.\ the proof of
Lemma~\ref{eltstrivonWF}~(\romannumeral1).)

We conclude from this that, if $W_F(X)$ consists of Hodge classes, then
$$\begin{array}{c}
W_F(X)\ \mbox{consists of}\\
\mbox{decomposable Hodge classes}
\end{array}
\,\Longleftrightarrow\,
\begin{array}{c}
\mbox{each of the spaces}\ W_F(Y_i^{m_i})\ \mbox{consists}\\
\mbox{of decomposable Hodge classes}
\end{array}\, .$$
\vspace{\baselineskip}

\noindent\refstepcounter{section}%
\textbf{\thesection.~Example.}\quad Let $X=Y_1 \times Y_2$ where the ratios 
$2\dim(Y_1)/[F:\Q]$ 
and $2\dim(Y_2)/[F:\Q]$ are odd. Then non-zero elements of $W_F(Y_1)$ and 
$W_F(Y_2)$ are not Hodge classes and therefore all non-zero elements of $W_F(X)$ 
are exceptional Hodge classes. Shioda's example mentioned above is of this type. 
See also \cite{SZ2}.
\vspace{\baselineskip}

\noindent\refstepcounter{section}%
\textbf{\thesection.}\quad We may from now on restrict our attention to the case 
that $k = 1$.
Since everything only depends on $X$ up to isogeny, we may even assume that $X
= Y^m$ for some $m \geq 1$, where $Y$ is simple. Let $D = \Eo(Y)$, let $E$ be
the center of $D$, and let $E_0$ be the maximal totally real subfield of $E$.
We write $e_0 = [E_0:\Q]$, $e=[E:\Q]$, $d^2 = \dim_E(D)$, and we say that $X$
and $Y$ are of type \hrom{1}, \hrom{2}, \hrom{3} or
\hrom{4} if the algebra $D$ is of the corresponding type in the Albert
classification.
\vspace{\baselineskip}

\noindent\refstepcounter{section}%
\textbf{\thesection.}\quad Choose a polarization $\lambda$ of $X$, and write 
$\alpha \mapsto \alpha^{\dagger}$ for the associated Rosati involution of 
$\Eo(X) \cong \Matr_m(D)$. Let $\cS_\lambda \subseteq \Eo(X)$ be the set of
$\dagger$-symmetric elements. We define the algebra $B \subseteq \Eo(X)$ as the
$\Q$-subalgebra generated by $\cS_\lambda$. If $\dagger^\prime$ is the Rosati
involution associated to another polarization, then $\dagger^\prime$ is
conjugated to $\dagger$ by an element of $\cS_\lambda$. It follows that the
algebra $B$ does not depend on the choice of the polarization $\lambda$.

For all possible types in the Albert classification one can determine the
algebra $B$ and its center $K_B$. The results are listed in Table~\ref{table1}.

\begin{table}[htbp]
%\label{table1}
\begin{center}
\begin{tabular}{|c|c|c|}
\hline
{\sl Type}\rule[-6pt]{0mm}{20pt} & $m=1$ & $m\geq2$\\
\hline
\hrom{1}\rule[-6pt]{0mm}{20pt} & \multicolumn{2}{c|}{$K_B = E$,\quad $B =
\Matr_m(E)$}\\
\hline
\hrom{2}\rule[-6pt]{0mm}{20pt} & \multicolumn{2}{c|}{$K_B = E$,\quad $B =
\Matr_m(D)$}\\
\hline
\hrom{3}\rule[-6pt]{0mm}{20pt} & $K_B = B = E$ & $K_B = E$,\quad $B =
\Matr_m(D)$\\
\hline
\hrom{4}, $d=1$\rule[-6pt]{0mm}{20pt} & $K_B = B = E_0$ & $K_B = E$,\quad $B =
\Matr_m(E)$\\
\hline
\hrom{4}, $d\geq2$\rule[-6pt]{0mm}{20pt} & \multicolumn{2}{c|}{$K_B = E$,\quad
$B = \Matr_m(D)$}\\
\hline
\end{tabular}

\end{center}
\caption{\label{table1}}
\end{table}

Let $\varphi = \varphi_X \colon V \times V \rightarrow \Q$ be the nondegenerate
alternating bilinear form associated to $\lambda$. We define the algebraic
group $\Gdiv(X) \subseteq \SP(V,\varphi)$ as the centralizer of $B$ in
$\SP(V,\varphi)$. More precisely,
$$\Gdiv(X) := \Gl_B(V) \cap \SP(V,\varphi)\, .$$
Of course, the main motivation for introducing this group $\Gdiv(X)$ is the
fact that it is the largest algebraic subgroup of $\Gl(V)$ defined over $\Q$
which leaves invariant all divisor classes in $\HH^2(X,\Q) = \bigwedge^2 V$. In
fact, the divisor classes, viewed as alternating bilinear forms on $V$, are
precisely the forms
$$\varphi_s \colon (v_1,v_2) \mapsto \varphi(s \cdot v_1, v_2)$$
for $s \in \cS_\lambda$.

The group $\Gdiv(X)$ does not depend on the choice of $\lambda$. Without loss
of generality we may therefore assume that there is a polarization $\mu$ on $Y$
such that $\lambda$ is the product polarization $\mu^m$ on $X=Y^m$. We write
$\varphi_Y$ for the alternating form on $V_Y$ associated to $\mu$ and $d
\mapsto d^\ast$ for the Rosati involution on $D=\Eo(Y)$. (With $\Eo(X) =
\Matr_m(D)$, the Rosati involutions $\ast$ and $\dagger$ are related by
$\alpha^\dagger = (\alpha_{ji}^\ast)$ for $\alpha = (\alpha_{ij}) \in
\Matr_m(D)$.)

{}From Table~\ref{table1} we see that in all cases there is a simple
$\Q$-subalgebra $\Delta \subseteq D = \Eo(Y)$ such that $\ast$ restricts to a
positive involution on $\Delta$, and such that $\Gdiv(X)$ is the centralizer of
$\Delta$ in $\SP(V_Y,\varphi_Y)$, embedded diagonally into $\SP(V,\varphi)$.
(In fact, if $m\geq 2$ or if $X$ is of type \hrom{1} or \hrom{2}, then we
simply have $\Delta = D$. If $m=1$ then $\Delta = B$.) This means that
$\Gdiv(X)$ is an algebraic group of the type that is usually studied in the
context of moduli problems of PEL-type. We refer the reader to \cite{Kott}, 
\cite{Haz1} or \cite{Kum1} for more information; here we only give a 
description of $\Gdiv(X) \otimes \C$.

Write $\Sigma_{E_0}$ for the set of complex (in fact real) embeddings of $E_0$.
We have
$$\SP(V_Y,\varphi_Y) \otimes \C = \prod_{\tau \in \Sigma_{E_0}}
\SP(V_{Y,\C}^{(\tau)}, \varphi_Y^{(\tau)})$$
and
$$\Delta \otimes \C =  \prod_{\tau \in \Sigma_{E_0}} \Delta_\C^{(\tau)}\, ,$$
where $\Delta_\C^{(\tau)}$ is a semi-simple $\C$-subalgebra of
$\End\big(V_{Y,\C}^{(\tau)}\big)$. We thus see that $\Gdiv(X) \otimes \C$
splits as the direct product of $e_0$ factors $\Gdiv^{(\tau)}$. Writing $k =
\dim_\Delta(V_Y)$, we have the following description of these factors and their
representations $V_{Y,\C}^{(\tau)}$. In each case $\Stand$ denotes the standard
representation of the group in question.

\begin{table}[htbp]
\begin{center}
\begin{tabular}{|c|c|ll|}
\hline
{\sl Type}\rule[-6pt]{0mm}{20pt} & $k$ &
\multicolumn{2}{c|}{$\Gdiv^{(\tau)}\quad\mbox{and}\quad
V_{Y,\C}^{(\tau)}$\quad}\\
\hline
\hrom{1}\rule[-6pt]{0mm}{20pt} & $2g/me$ & $\Gdiv^{(\tau)} \cong \SP_{2k}$ &
$V_{Y,\C}^{(\tau)}\cong \Stand$\\
\hline
\hrom{2}\rule[-6pt]{0mm}{20pt} & $g/2me$ & $\Gdiv^{(\tau)} \cong \SP_{2k}$ &
$V_{Y,\C}^{(\tau)}\cong \Stand \oplus \Stand$\\
\hline
\hrom{3}, $m=1$\rule[-6pt]{0mm}{20pt} & $2g/e$ & $\Gdiv^{(\tau)} \cong
\SP_{2k}$ & $V_{Y,\C}^{(\tau)}\cong \Stand$\\
\hline
\hrom{3}, $m\geq 2$\rule[-6pt]{0mm}{20pt} & $g/2me$ & $\Gdiv^{(\tau)} \cong
\OO_{2k}$ & $V_{Y,\C}^{(\tau)}\cong \Stand \oplus \Stand$\\
\hline
\hrom{4}, $d=1$, $m=1$\rule[-6pt]{0mm}{20pt} & $2g/e_0$ & $\Gdiv^{(\tau)} \cong
\SP_{2k}$ & $V_{Y,\C}^{(\tau)}\cong \Stand$\\
\hline
\hrom{4}, $d\geq2$ or $m \geq 2$\rule[-6pt]{0mm}{20pt} & $2g/med^2$ &
$\Gdiv^{(\tau)} \cong \Gl_{dk}$ & $V_{Y,\C}^{(\tau)}\cong \Stand \oplus
\Stand^\vee$\\
\hline
\end{tabular}

\end{center}
\caption{\label{table2}}
\end{table}
\vspace{\baselineskip}

\noindent\refstepcounter{section}\label{Gdivprops}%
\textbf{\thesection.~Lemma.}\quad (\romannumeral1) \textsl{The center of $\Gdiv(X)$ 
is the group $\Uu_{K_B}$ given by
$$\Uu_{K_B}(R) = \{a \in (K_B \otimes_\Q R)^\ast \mid a a ^\dagger = 1\}\, .$$
For $X$ of type \hrom{4} with either $d \geq 2$ or $m \geq 2$ this is a
connected torus of rank $e_0$; in all other cases it is finite.}

(\romannumeral2) \textsl{If $X$ is not of type \hrom{3} with $m \geq 2$, then
$\Gdiv(X)$ is geometrically connected; for $X$ of type \hrom{3} and $m \geq 2$,
the group $\pi_0\big(\Gdiv(X)\big)$ has (geometrically) order $2^{e_0}$.}

(\romannumeral3) \textsl{$\End(V_X)^{\Gdiv(X)} = B$; $(\bigwedge^2 V_X)^{\Gdiv(X)} =
\cB^1(X)$, and $(\oplus_i \bigwedge^i V_X)^{\Gdiv(X)} = \calD^\punt(X)$.}
\vspace{\baselineskip}

\noindent
{\bf Proof.} To prove this, we can first extend scalars to $\C$, and
(\romannumeral1) and (\romannumeral2) then readily follow from
Table~\ref{table2}. The last statement is based on results from classical
invariant theory; both the statement and its proof are in fact easy variants of
\cite[Lemma~3.6]{Kum1} (see also \cite{Haz1}). \Qed
\vspace{\baselineskip}

\noindent\refstepcounter{section}\label{eltstrivonWF}%
\textbf{\thesection.~Lemma.}\quad (\romannumeral1) \textsl{The space 
$W_F = \bigwedge^r_F V$ can
naturally be identified with a subspace of $\bigwedge^r_\Q V$.}

(\romannumeral2) \textsl{If $g \in \Gl_\Q(V)$ acts as the identity on $W_F$ then 
$g$ is $F$-linear, hence $g \in \Sl_F(V)$.}
\vspace{\baselineskip}

\noindent
{\bf Proof.} There is a canonical isomorphism $\Tr_{F/\Q} \colon \Hom_F(V_X,F)
\isomarrow \Hom_\Q(V_X,\Q)$ and we simply write $V_X^\vee$ for this space.
There is a canonical surjective map $\bigwedge^r_\Q V_X^\vee \twoheadrightarrow
\bigwedge^r_F V_X^\vee$. It induces an injective linear map
$$\Alt_F^r(V_X) \cong \Hom_F(\bigwedge_F^r V_X^\vee,F) \cong
\Hom_\Q(\bigwedge_F^r V_X^\vee,\Q) \hooklongrightarrow \Hom_\Q(\bigwedge_\Q^r
V_X^\vee,\Q) \cong \Alt_\Q^r(V_X)\, ,$$
where all isomorphisms are canonical. Since our fields are of characteristic 0,
there are identifications $\bigwedge_F^r V_X \cong \Alt_F^r(V_X)$ and
$\bigwedge_\Q^r V_X \cong \Alt_\Q^r(V_X)$, and statement (\romannumeral1)
follows.

(\romannumeral2) Choose an $F$-basis for $V_X$. This gives an isomorphism $F
\isomarrow \Hom_F(\bigwedge_F^r V_X^\vee, F)$ by sending $f \in F$ to the
functional $t_1 \wedge_F \ldots \wedge_F t_r \mapsto f \cdot
\determ_F(t_1,\ldots,t_r)$. Suppose $g \in \Gl_\Q(V_X)$ acts trivially on
$W_F$. This means that for all $f \in F$ and all $r$-tuples $t_1, \ldots, t_r
\in V^\vee$ we have the identity
$$\Tr_{F/\Q}\big(f \cdot \determ_F(t_1,\ldots,t_r)\big) = \Tr_{F/\Q}\big(f
\cdot \determ_F(g \cdot t_1,\ldots,g \cdot t_r)\big)\, .$$
The trace form being non-degenerate it follows that
$$\determ_F(t_1,\ldots,t_r) = \determ_F(g \cdot t_1,\ldots,g \cdot
t_r)\qquad\mbox{for all $t_1, \ldots, t_r \in V_X^\vee$}\, .$$
For $f \in F$ this gives the identity
$$\determ_F(g \cdot (f t_1),g \cdot t_2,\ldots,g \cdot t_r) = f \cdot
\determ_F(g \cdot t_1,\ldots,g \cdot t_r) = \determ_F(f (g \cdot t_1),g \cdot
t_2,\ldots,g \cdot t_r)\, ,$$
and (\romannumeral2) readily follows from this. \Qed
\vspace{\baselineskip}

\noindent\refstepcounter{section}\label{crit2}%
\textbf{\thesection.~Criterion.}\quad {\slshape Suppose $X$ is isogenous to a power 
$Y^m$ of a simple abelian variety $Y$ (see
section~\ref{decompuptoisog}).
Suppose $F \hookrightarrow \Eo(X)$ is a subfield such that $W_F = \bigwedge^r
V_X$ consists of Hodge classes (see section~\ref{crit1}, and recall that this assumption
implies that $r = \dim_F(V_X)$ is even). Then either all classes in $W_F$ are
decomposable, or all non-zero classes in $W_F$ are exceptional; this last
possibility occurs precisely in the following cases:

$Y$ is of Type \hrom{3}, $m=1$ and $F \subsetneq E$,

$Y$ is of Type \hrom{3}, $m \geq 2$ and the integer $2m \cdot [E:\Q] / [F:\Q]$
is odd,

$Y$ is of Type \hrom{4}, $d=1$, $m=1$ and $F \subsetneq E_0$,

$Y$ is of Type \hrom{4} with $d\geq2$ or $m \geq 2$ and}
\begin{equation} \mbox{\textsl{the map}}\quad \theta \colon E_- \hooklongrightarrow
\End_F(V_X) \stackrel{\Tr_F}{\verylongarrow} F\quad \mbox{\textsl{is non-zero}}\, .
\label{eq:Type4cond}\end{equation}
\vspace{\baselineskip}

\noindent
{\bf Proof.} First assume that $X$ is either of type \hrom{1}, \hrom{2} or of
type \hrom{3} with $m = 1$, or that $X$ is of type \hrom{4} with $d=1$ and
$m=1$. We claim that, in these cases, $\Gdiv(X)$ acts as the identity on $W_F$
if and only if $F \subseteq B$. In the ``only if'' direction, this follows from
Lemmas~\ref{eltstrivonWF}~(\romannumeral2) and
\ref{Gdivprops}~(\romannumeral3). Conversely, suppose that $F \subseteq B$, so
that $\Gdiv(X) \subseteq \Gl_F(V_X)$. In the cases we are considering, the
group $\Gdiv(X)$ is connected and semi-simple, so $\Gdiv(X) \subseteq
\Sl_F(V_X)$, hence $\Gdiv(X)$ acts trivially on $W_F$.

Next assume that $X$ is of type \hrom{4} with either $m \geq 2$ or $d \geq 2$.
We have $F \subseteq B = \Eo(X)$. Since in this case $\Gdiv(X)$ is connected
(see \ref{Gdivprops}), it acts trivially on $W_F$ if and only if the
composition
$$ \Uu_E = Z(\Gdiv) \subset \Gdiv(X) \hooklongrightarrow \Gl_F(V_Y)
\stackrel{\determ_F}{\verylongarrow} F^\ast$$
is trivial. The torus $\Uu_E$ being connected, this is the case if and only if
the induced map on Lie algebras
$$\theta \colon E_- \hooklongrightarrow \End_F(V_X)
\stackrel{\Tr_F}{\verylongarrow} F$$
is zero.

Finally we consider the case where $X$ is of type \hrom{3} with $m \geq 2$. We
have $B = \Matr_m(D)$, hence $F \subseteq B$. The centralizer $\cC$ of $B$ in
$\End_E(V_X)$ is isomorphic to $\Matr_k(D^\prime)$, where $D^\prime$ is the
opposite algebra of $D$ and $k = \dim_E(V_X) / (m \cdot \dim_{E}(D)) = g /
2me$. Let $\Nrd \colon \cC \rightarrow E$ denote the reduced norm.

Let us first assume that $F$ contains the field $E$. The restriction of the
$F$-linear determinant map $\determ_F \colon \End_F(V_X) \rightarrow F$ to the
subalgebra $\cC \subseteq \End_F(V_X)$ coincides with a power $\Nrd^q$ of the
reduced norm, and comparing the degrees we find that $q = \dim_F(V_X)/2k =
2m/[F:E]$.

Let us now show that, in this case, $W_F$ consists of decomposable classes if
and only if $2m/[F:E]$ is even. By the same arguments as in the type \hrom{1}
and \hrom{2} case, we see that the connected component of the identity
$\Gdiv^0$ acts trivially on $W_F$. Next consider an embedding $\tau \colon E
\rightarrow \C$, and an element $g \in \Gdiv^{(\tau)}(\C)$ which is not in the
connected component $(\Gdiv^{(\tau)})^0$. The $\C$-linear extension of the
reduced norm to $\cC \otimes_{E,\tau} \C \cong \Matr_{2k}(\C)$ is nothing but
the $\C$-linear determinant, so $\Nrd_\C(g) = \determ_\C(g) = -1$. (Recall that
$\Gdiv^{(\tau)} \cong \OO_{2k}$.) It then follows from the preceding remarks
that $g$ acts on $W_F \otimes \C$ as multiplication by $(-1)^{2m/[F:E]}$, which proves the
assertion.

Next let us do the general case (i.e., no longer assuming that $E \subseteq
F$). The compositum $EF \subseteq \Eo(X)$ is a product of fields, say $EF = K_1
\times \cdots \times K_t$. Correspondingly, we can decompose $V_X$ as a direct
sum $V_1 \oplus \cdots \oplus V_t$, and we have
$$W_F = \determ_F(V_X) = \determ_F(V_1) \otimes_F \cdots \otimes_F
\determ_F(V_t)\, ,$$
where $F$ acts on $V_i$ through its embedding into $K_i$. For $g \in \Gdiv$ we
have
$$\determ_F(g;V_i) = \Norm_{K_i/F} \big(\determ_{K_i}(g;V_i)\big)\, .$$
To calculate $\determ_{K_i}(g;V_i)$ we can argue as in the preceding case
(replacing $F$ by $K_i$ in the argument). If $g$ does not lie in the connected
component of $\Gdiv$ we find that 
$\determ_{K_i \otimes \C}(g;V_i \otimes \C) = (-1)^{q_i}$, where
$q_i = 2m \cdot \dim_{K_i}(V_i) / \dim_E(V_X)$. In total we therefore get
$$\determ_F(g;V_X \otimes \C) = \prod_{i=1}^t (-1)^{[K_i:\Q]\cdot q_i} =
(-1)^{\sum_{i=1}^t 2m \cdot [K_i:F] \cdot \dim_{K_i}(V_i) / \dim_E(V_X)}\, ,$$
and since
$$\sum_{i=1}^t \frac{2m \cdot [K_i:F] \cdot \dim_{K_i}(V_i)}{\dim_E(V_X)} =
\frac{2m}{\dim_E(V_X)} \cdot \sum_{i=1}^t \dim_F(V_i) = 2m \cdot
\frac{\dim_F(V_X)}{\dim_E(V_X)} = 2m \cdot \frac{[E:\Q]}{[F:\Q]}\, ,$$
we see that $W_F$ consists of decomposable Hodge classes if and only if the
integer $2m [E:\Q]/[F:\Q]$ is even. This proves the criterion. \Qed
\vspace{\baselineskip}

\noindent\refstepcounter{section}%
\textbf{\thesection.~Remark.}\quad The existence of exceptional Weil classes 
on simple abelian varieties of type \hrom{3} was proven by V. K. Murty \cite{Kum1}.
\vspace{\baselineskip}

\noindent\refstepcounter{section}%
\textbf{\thesection.}\quad In practice, the condition (\ref{eq:Type4cond}) 
in the case that $Y$ is of
Type~\hrom{4} with $d \geq 2$ or $m \geq 2$, is not as easy to verify as the
conditions on $F$ in the other cases. Let us add some remarks to illustrate
what can happen. Here, as before, we assume that $X$ is isogenous to $Y^m$,
where $Y$ is simple and $m \geq 1$.

(\romannumeral1)~Suppose we have two subfields $F \subseteq F^\prime \subseteq
\Eo(X)$. As remarked before, if $X$ is of type \hrom{1}, \hrom{2} or \hrom{3},
then the spaces $W_F$ and $W_{F^\prime}$ always consist of Hodge classes. In
the type \hrom{4} case this is not true in general. However, it follows
directly from Criterion~\ref{crit1} that
$$W_{F^\prime}\ \mbox{consists of Hodge classes}\quad \Longrightarrow\quad 
W_F\ \mbox{consists of Hodge classes}\, .$$
Moreover, if $F^\prime \supseteq E$, then the converse is true. To see this,
let us recall that the center $Z(\Hdg)$ of the Hodge group is contained in the
torus $\Uu_E$, and that the action of $\Hdg$ on $W_F$ is given by the
$F$-linear determinant. Therefore, if $E \subset F^\prime$, then $W_{F^\prime}$
consists of Hodge classes if and only if $\Hdg$ is semi-simple. If this holds,
then for any other subfield $F \subseteq \Eo(X)$, the space $W_F$ also consists
of Hodge classes.

Furthermore, if $F \subseteq F^\prime$ then we have the implication
$$
\begin{array}{c}
W_{F^\prime}\ \mbox{consists of}\\
\mbox{decomposable Hodge classes}
\end{array}
\quad\Longrightarrow\quad
\begin{array}{c}
W_F\ \mbox{consists of}\\
\mbox{decomposable Hodge classes}
\end{array}\, .$$
This is a direct consequence of Criterion~\ref{crit2}. It can also be seen more
directly, by using that $W_F$ is contained in the vector subspace generated by
exterior products of elements of $W_{F^\prime}$.

(\romannumeral2)~Assume $X$ is of type~\hrom{4} with either $d \geq 2$ or $m
\geq 2$. If the map $\theta \colon E_- \hooklongrightarrow \End_F(V_X)
\stackrel{\Tr_F}{\verylongarrow} F$ is zero, then the intersection $E \cap F$
is a totally real subfield of $E$, i.e., $E \cap F \subseteq E_0$. In fact, if
$E \cap F$ is not totally real, then (being a CM-subfield of $E$) it must
contain totally imaginary elements $0 \neq \alpha \in E_-$, which then
obviously have a non-zero trace over $F$.

(\romannumeral3)~Let us show that the converse of (\romannumeral2) holds if
either $E$ or $F$ is Galois over $E \cap F$. So, we assume that $E \cap F$ is
totally real and that $E$ is Galois over $E \cap F$. The compositum $EF
\subseteq \Eo(X)$ is a product of fields, say $EF = K_1 \times \cdots \times
K_t$. Correspondingly, we can decompose $\HH^1(X,\Q)$ as a direct sum $V_1
\oplus \cdots \oplus V_t$. It suffices to show that each of the maps
$\theta_i \colon E_- \hooklongrightarrow \End_F(V_i)
\stackrel{\Tr_F}{\verylongarrow} F$ is zero.

We have $E_- \hooklongrightarrow K_i \subset \End_{K_i}(V_i)$, and it follows
that for $\alpha \in E_-$,
$$\theta_i(\alpha) = \Tr_{K_i/F}(\dim_{K_i}(V_i) \cdot \alpha) =
\dim_{K_i}(V_i) \cdot \Tr_{K_i/F}(\alpha)\, .$$
On the other hand, if either $E$ or $F$ is Galois over $E \cap F$, then the
trace of $\alpha$ (considered as an element of $K_i$) over $F$ lies in $E \cap
F$, hence
$$\Tr_{K_i/F} (\alpha) = \frac{[K_i:E]}{[F: E\cap F]} \cdot \Tr_{E/E\cap
F}(\alpha) = 0\, .$$

(\romannumeral4)~By means of an example, we can show that in general the
converse of (\romannumeral2) does not hold. For this, suppose that we have a
field $K$, containing a CM-subfield $E$ and a subfield $F$ such that $E \cap F$
is totally real and such that the map $\theta \colon E_- \subset K
\stackrel{\Tr_F}{\verylongarrow} F$ is non-zero. Then we obtain an example of
the kind we are looking for by using the constructions and results of
\cite{Shim}. First we choose a central simple algebra $D$ over $E$ containing
$K$ as a subfield, and then we take an abelian variety $X$ with $\Eo(X) = D$
and such that $n_\sigma = n_{\sigma^\prime}$ for all $\sigma \in \Sigma_E$
(which is possible, see \cite[Thm.~5]{Shim}). Notice that this $X$ is of
type~\hrom{4} with $m = 1$ and $d \geq 2$ and that $W_F$ consists of Hodge
classes (using remark (\romannumeral1) above).

The construction of the fields $K$, $E$ and $F$ can be done using Galois
theory. For example, we can start with a CM-field $K$ which is Galois over $\Q$
with group $\{\pm 1\}^n \rtimes S_n$ (complex conjugation corresponding to
$(-1,\ldots,-1) \rtimes \Id$), take $E$ to be the CM-subfield of elements that
are invariant under $\{\pm 1\}^{n-1} \rtimes S_{n-1}$, and let $F$ be the fixed
field of the transposition $(1,\ldots,1) \rtimes (n-1\quad n)$.

(\romannumeral5)~In the opposite direction, we can also use the results of
\cite{Shim} to construct examples where the space $W_F$ consists of exceptional
Hodge classes. In particular, we see that whenever we have number fields $E
\subset K \supset F$ with $E$ a CM-field and $E \cap F$ not totally real, then
there exists a simple abelian variety $X$ of type \hrom{4} with $d \geq 2$,
such that $F \hookrightarrow \Eo(X)$ and such that the associated space $W_F$
consists of exceptional Hodge classes.
\vspace{\baselineskip}

\noindent\refstepcounter{section}\label{Q1Q2Tate}%
\textbf{\thesection.}\quad To conclude, let us consider the questions Q1 and Q2 
(see the introduction) in the context of Tate
classes. So, let $X$ be an abelian variety defined over a number field $K$,
assumed to be large enough, as always. Write $V_\ell = V_{\ell,X}$ for the
first \'etale cohomology $\HH^1_{\et}(X_{\Qbar},\Ql)$. If $F$ is a subfield of
$\Eo(X)$, then the space
$$W_{\ell,F} = \bigwedge^r_{F\otimes \Ql} V_{\ell,X}$$
can be identified with a subspace of $\HH^r_{\et}(X_{\Qbar},\Ql)$.

We claim that the criteria \ref{crit1} and \ref{crit2} (replacing ``Hodge
classes'' by ``Tate classes'' and $W_F$ by $W_{F,\ell}$) are valid in this
context as well. The proof of Criterion \ref{crit1} is essentially the same;
the main difference is that one uses a Hodge-Tate decomposition to replace the
Hodge decomposition (cf.\ the proof of \cite[Lemma~2.8]{MZ}).

Criterion~\ref{crit2} then follows from the fact that the rings of divisor
classes are ``the same'' in the ``Hodge'' and the ``Tate'' context, i.e.,
$(\calD^\punt(X_\C)) \otimes \Ql = \calD_\ell^\punt(X_{\Qbar})$. This, of course,
is an immediate consequence of the fact that $(\cB^1(X_\C)) \otimes \Ql =
\cB_\ell^1(X_{\Qbar})$ (by the Lefschetz theorem on (1,1) classes and a special
case of the Tate conjecture, as proven by Faltings in \cite{Falt}).

\vspace{2\baselineskip}

\textsc{\footnotesize B.J.J.~Moonen, Westf\"alische Wilhelms-Universit\"at
M\"unster, Mathematisches
Institut, Einsteinstra\ss e 62, 48149 M\"unster, Germany.}

\textsl{\footnotesize Email:\/}  \texttt{\footnotesize moonen@math.uni-muenster.de}
\vspace{.5\baselineskip}

\textsc{\footnotesize Yu.G.~Zarhin, Department of Mathematics, Pennsylvania State University, 
University Park, PA 16802, USA}

\textsc{\footnotesize Institute for Mathematical Problems in Biology, 
Russian Academy of Sciences, Pushchino, Moscow Region, 142292, Russia}

\textsl{\footnotesize Email:\/}  \texttt{\footnotesize zarhin@math.psu.edu}
\end{document}